\documentclass[12pt]{article}

\usepackage{amsmath}\usepackage{amssymb}\usepackage{amsfonts}\usepackage{amsthm}
\usepackage{color}
\usepackage{graphicx}
\usepackage[english]{babel}\usepackage[latin1]{inputenc}\usepackage{times}\usepackage[T1]{fontenc}

\newtheorem{thm}{\bf{Theorem}}[section]
\newtheorem{lem}[thm]{\bf{Lemma}}

\newtheorem{df}[thm]{\bf{Definition}}
\newtheorem{cor}[thm]{\bf{Corollary}}

\newtheorem{prop}[thm]{\bf{Proposition}}

\newtheorem{ex}[thm]{\bf{Example}}

\numberwithin{equation}{section}

\newcommand{\dom}{\operatorname{dom}}
\newcommand{\gph}{\operatorname{gph}}
\newcommand{\dist}{\operatorname{dist}}
\newcommand{\conv}{\operatornamewithlimits{conv}}
\newcommand{\argmin}{\operatornamewithlimits{argmin}}

\title{Parametrically Prox-Regular Functions\\}
\date{\today}
\author
{W. Hare\thanks{University of British Columbia, Okanagan Campus (UBCO), 3333 University Way, Kelowna, BC, Canada.   Research by this author was supported by NSERC of Canada. {\tt warren.hare@ubc.ca}} \and C. Planiden\thanks{UBCO. Research by this author was supported by UBC UGF and by NSERC of Canada. {\tt chayneplaniden@hotmail.com}} }

\begin{document}

\maketitle\author
\setcounter{page}{1}\pagenumbering{arabic}
\begin{abstract}
Prox-regularity is a generalization of convexity that includes all $\mathcal{C}^2$, lower-$\mathcal{C}^2$, strongly amenable, and primal-lower-nice functions. The study of prox-regular functions provides insight on a broad spectrum of important functions. Parametrically prox-regular (para-prox-regular) functions are a further extension of this family, produced by adding a parameter. Such functions have been shown to  play a key role in understanding stability of minimizers in optimization problems. This document discusses para-prox-regular functions in $\mathbb{R}^n.$\\
 \indent We begin with some basic examples of para-prox-regular functions, and move on to the more complex examples of the convex and nonconvex proximal average. We develop an alternate representation of a para-prox-regular function, related to the monotonicity of an $f$-attentive $\epsilon$-localization as was done for prox-regular functions \cite{proxfunc}. This extends a result of Levy \cite{calmmin}, who used an alternate approach to show one implication of the relationship (we provide a characterization). We analyze two common forms of parametrized functions that appear in optimization: finite parametrized sum of functions, and finite parametrized max of functions. The example of strongly amenable functions by Poliquin and Rockafellar \cite{calcprox} is given, and a relaxation of its necessary conditions is presented.
\end{abstract}

\section{Introduction}
\par Convex analysis has long been a prominent area of study in mathematics, and in optimization in particular. One method of unifying several classes of functions together with convex functions, that allows us to take advantage of their different properties, is through the notion of prox-regular functions.
\par Prox-regularity is a local property defined as follows. A function $f$ is said to be prox-regular at $\bar{x}$ for a subgradient $\bar{v}$ if there exist $\epsilon>0$ and $r\geq0$ such that
$$f(x')\geq f(x)+\langle v,x'-x\rangle-\frac{r}{2}|x'-x|^2$$
under $\epsilon$-neighborhoods of $\bar{x}$, $f(\bar{x})$ and $\bar{v}$ (a formal definition appears in Section \ref{section:prelim}). Prox-regularity was first introduced in the mid-nineties \cite{proxsmooth,genhess,proxfunc}, when it was shown that all proper convex functions are prox-regular, as are lower-$\mathcal{C}^2$, strongly amenable, and primal-lower-nice functions \cite{proxfunc}. The class of lower-$\mathcal{C}^2$ functions includes the class of $\mathcal{C}^2$ functions, thus the study of prox-regular functions gives us a way to connect the analyses of smooth functions and convex functions. It is quite common to see results that begin with a function of one of the classes mentioned above, and provide a related function that is of another of those classes. For example, Poliquin and Rockafellar \cite{proxfunc} showed that if a function is prox-regular, then its Moreau envelope is a $\mathcal{C}^{1+}$ function.
\par Much of the theory in the field of prox-regularity has been extended from finite-dimensional space to Hilbert and Banach spaces \cite{hilbert,banach}. Prox-regularity is key in matters of function stability \cite{tiltstab,subdiffstab,proxmap,tiltopt}, and it has been shown to extend results previously known for convex functions, such as the convergence rates of variational inequalities \cite{varineq}, alternating and averaged projection algorithms \cite{aveproj}, active manifold identification \cite{proxmeth}, proximal averages \cite{proxave}, and proximal bundle methods \cite{proxbund}.
\par Parametric prox-regularity, first introduced in 2000 by Levy, Poliquin and Rockafellar \cite{stability}, extends prox-regularity to parametrized functions. The formal definition of parametric prox-regularity is presented in Definition \ref{df:para}. Little research has been done on parametrically prox-regular (henceforth para-prox-regular) functions to date; the purpose of this article is to take a step in that direction. We present an alternate definition similar to that of prox-regular functions given by Poliquin and Rockafellar \cite{proxfunc}. We examine several common styles for constructing parametrized functions, and explore where the results are para-prox-regular.
\par The organization of this article is as follows. Section \ref{section:prelim} starts with the definitions of subgradients, prox-regular functions, and para-prox-regular functions. It also provides some preliminary results for use in subsequent sections. Section \ref{section:ex} presents some examples of prox-regular and para-prox-regular functions. These include the proximal average \cite{proxpoint,howto} and the NC-proximal average \cite{proxave}, which are parametrized functions designed to transform one function into another in a continuous manner. Section \ref{section:mono} contains an alternate definition of para-prox-regular functions. This result is similar to Poliquin's and Rockafellar's alternate definition of prox-regular functions, which is based on an $f$-attentive $\epsilon$-localization $T$ of $\partial f$ around $(\bar{x},\bar{v})$ (the prox-regularity of $f$ is equivalent to the monotonicity of $T+rI$). Section \ref{section:main} presents three methods of constructing parametrized functions, and discusses when each method results in a para-prox-regular function. As a corollary, we answer an open question posed in \cite{proxave}. The methods presented are scalar multiplication of a function by a positive parameter, the parametrized sum of a finite collection of functions, and the finite parametrized max of continuous functions. An expansion of \cite[Theorem 3.2]{calcprox} is also seen in Section \ref{section:main}, extending the required two initial prox-regular functions to any finite number of prox-regular functions. Section \ref{sec:open} poses some suggestions for future research.

\section{Preliminaries}
\label{section:prelim}
This section includes the definitions of prox-regular and para-prox-regular functions. Henceforth, functions are assumed to be proper and lower semicontinuous (lsc) as defined in \cite{rockwets}, unless otherwise stated.
\begin{df}{(Little-o)} For two real functions $f$ and $g$, we say $f(x)=o(g(x))$ at $x_0$ if and only if $\lim\limits_{x\rightarrow x_0}\frac{f(x)}{g(x)}=0$.
\end{df}
The first important definition we will need is that of \emph{subgradients}. There are several types of subgradients; their definitions are restated here.
\begin{df}\cite[Def. 8.3]{rockwets}
\label{df:subgrad}
Let $f:\mathbb{R}^n\rightarrow\mathbb{R}\cup\{\infty\}$, $\bar{x}\in\dom f$. A vector $\bar{v}\in\mathbb{R}^n$ is
\begin{itemize}
\item[a)] a \textbf{regular subgradient} of $f$ at $\bar{x}$, written $\bar{v}\in\hat{\partial}f(\bar{x})$, if
\begin{equation}
f(x)\geq f(\bar{x})+\langle\bar{v},x-\bar{x}\rangle+o(|x-\bar{x}|).
\end{equation}
\item[b)] a \textbf{(general) subgradient} of $f$ at $\bar{x}$, written $\bar{v}\in\partial f(\bar{x})$, if there exist sequences $\{x^\nu\}_{\nu=1}^\infty\subseteq\dom f$ and $\{v^\nu\}_{\nu=1}^\infty\subseteq\hat{\partial}f(x^\nu)$ such that $(x^\nu,f(x^\nu))\rightarrow(\bar{x},f(\bar{x}))$ and $v^\nu\rightarrow\bar{v}$.
\item[c)] a \textbf{horizon subgradient} of $f$ at $\bar{x}$, written $v\in\partial^\infty f(\bar{x})$, if the same holds as in b) except that instead of $v^\nu\rightarrow v$ one has $\lambda^\nu v^\nu\rightarrow v$ for some sequence $\lambda^\nu\searrow0$, or in other words, $v^\nu\rightarrow dir(v)$ (or $v=0$).
\item[d)] a \textbf{proximal subgradient} of $f$ at $\bar{x}$ if it is a regular subgradient whose error term $o(|x-\bar{x}|)$ is $\frac{r}{2}|x-\bar{x}|^2$ for some $r>0$.
\end{itemize}
\end{df}
The set of all regular subgradients of $f$ at $\bar{x}$ is denoted $\hat{\partial}f(\bar{x})$ and is called the \emph{regular subdifferential} of $f$ at $\bar{x}$. Similarly, the set of (general) subgradients is denoted $\partial f(\bar{x})$ and is called the \emph{subdifferential} and the set of horizon subgradients is denoted $\partial^{\infty}f(\bar{x})$ and called the \emph{horizon subdifferential}. For a function $f:\mathbb{R}^m\times\mathbb{R}^n \rightarrow\mathbb{R}\cup\{\infty\}$, $f(x,\lambda)$, we shall also use $\partial_xf(\cdot,\cdot)$ to denote the subdifferential with respect to $x$.

Now we are ready for the definition of prox-regular functions.

\begin{df}\cite[Def. 13.27]{rockwets}
\label{df:pr1}
A proper function $f:\mathbb{R}^n\rightarrow\mathbb{R}\cup\{\infty\}$ that is finite at $\bar{x}$ is \textbf{prox-regular} at $\bar{x}$ for $\bar{v}$ with parameters $\epsilon>0, r\geq0$, where $\bar{v}\in\partial f(\bar{x})$, if and only if $f$ is locally lsc at $\bar{x}$ and
\begin{equation}
\label{eq:pr1}
f(x')\geq f(x)+\langle v,x'-x\rangle-\frac{r}{2}|x'-x|^2
\end{equation}
whenever $|x'-\bar{x}|<\epsilon$, $|x-\bar{x}|<\epsilon$, $|f(x)-f(\bar{x})|<\epsilon$, $v\in\partial f(x)$ and $|v-\bar{v}|<\epsilon$. If this holds true for all $\bar{v}\in\partial f(\bar{x})$, then $f$ is said to be prox-regular at $\bar{x}$ (without reference to any subgradient). If $f$ is prox-regular for all $\bar{x}\in \dom f$, then it is said to be a prox-regular function (without reference to any point).
\end{df}
By choosing $\epsilon$ slightly smaller, we may change to a strict inequality and rewrite as follows.

\begin{cor}\cite[Def. 1.1]{calcprox}
A proper function $f:\mathbb{R}^n\rightarrow\mathbb{R}\cup\{\infty\}$ is prox-regular at $\bar{x}$ for $\bar{v}$ with parameters $\epsilon>0$, $r>0$, if and only if $f$ is locally lsc at $\bar{x}$, $\bar{v}\in\partial f(\bar{x})$, and
\begin{equation}
\label{eq:pr2}
f(x')>f(x)+\langle v,x'-x\rangle-\frac{r}{2}|x'-x|^2
\end{equation}
whenever $x'\neq x$, $|x'-\bar{x}|<\epsilon$, $|x-\bar{x}|<\epsilon$, $|f(x)-f(\bar{x}|<\epsilon$, $|v-\bar{v}|<\epsilon$ and $v\in\partial f(x)$.
\end{cor}
Notice that since we now have a strict inequality for $r$, we must include the condition $x'\neq x$ to remove the possibility of $r=0$. Graphically speaking, due to the last term on the right-hand side of (\ref{eq:pr1}), we can think of prox-regular functions as those functions which are locally bounded below by a tangent concave quadratic function. An alternate definition of prox-regularity, which will be used in Section \ref{section:main}, is described in the following theorem.

\begin{thm}\cite[Theorem 3.2]{proxfunc}
\label{thm:prox1}
A proper function $f:\mathbb{R}^n\rightarrow\mathbb{R}\cup\{\infty\}$ is prox-regular at $\bar{x}$ for $\bar{v}$ with parameters $\epsilon>0, r\geq0$ if and only if $f$ is locally lsc at $\bar{x}$, $\bar{v}\in\partial f(\bar{x})$ and
$$\langle v'-v,x'-x\rangle\geq-r|x'-x|^2$$
whenever $|x'-\bar{x}|<\epsilon$, $|x-\bar{x}|<\epsilon$, $|f(x')-f(\bar{x})|<\epsilon$, $|f(x)-f(\bar{x})|<\epsilon$, $v'\in\partial f(x')$, $v\in\partial f(x)$, $|v'-\bar{v}|<\epsilon$, and $|v-\bar{v}|<\epsilon$.
\end{thm}
\medskip
From Definition \ref{df:pr1} we see that a prox-regular function is one that is locally bounded below by a quadratic function. Theorem \ref{thm:prox1} gives us an alternate expression that defines prox-regularity in terms of local pre-monotonicity of the subdifferential operator. This form is the inspiration for our main result in Section \ref{section:mono}.
\par Para-prox-regularity is a natural extension of prox-regularity that allows for a parametrized function.
\begin{df}\cite[Definition 2.1]{stability}
\label{df:para}
The proper function $f:\mathbb{R}^n\times\mathbb{R}^m\rightarrow\mathbb{R}\cup\{\infty\}$ is para-prox-regular in $x$ at $\bar{x}$ for $\bar{v}$ with compatible parametrization by $\lambda$ at $\bar{\lambda}$ (also said to be para-prox-regular in $x$ at $(\bar{x},\bar{\lambda})$ for $\bar{v}$), with parameters $\epsilon>0, r\geq0$, if and only if $f$ is locally lsc at $(\bar{x},\bar{\lambda})$, $\bar{v}\in\partial_xf(\bar{x},\bar{\lambda})$ and
\begin{equation}
\label{eq:para}
f(x',\lambda)\geq f(x,\lambda)+\langle v,x'-x\rangle-\frac{r}{2}|x'-x|^2
\end{equation}
whenever $|x'-\bar{x}|<\epsilon$, $|x-\bar{x}|<\epsilon$, $|f(x,\lambda)-f(\bar{x},\bar{\lambda})|<\epsilon$, $v\in\partial_xf(x,\lambda)$, $|v-\bar{v}|<\epsilon$, and $|\lambda-\bar{\lambda}|<\epsilon$.  The function $f$ is continuously para-prox-regular in $x$ at $(\bar{x},\bar{\lambda})$ for $\bar{v}$ if, in addition, $f(x,\lambda)$ is continuous as a function of $(x,v,\lambda)\in\gph\partial_xf$ at $(\bar{x},\bar{v},\bar{\lambda})$. It is para-prox-regular in $x$ at $(\bar{x},\bar{\lambda})$ (with no mention of $\bar{v}$) if it is para-prox-regular in $x$ at $(\bar{x},\bar{\lambda})$ for all $\bar{v}\in\partial_xf(\bar{x},\bar{\lambda})$, and it is a para-prox-regular function in $x$ (with no mention of $\bar{x}$ or $\bar{\lambda}$) if it is para-prox-regular in $x$ at $(\bar{x},\bar{\lambda})$ for all $(\bar{x},\bar{\lambda})\in\dom f$.
\end{df}
\medskip
A para-prox-regular function $f$ of $(x,\lambda)$ is a function that is prox-regular as a function of $x$ at $(\bar{x},\bar{\lambda})$, and continues to be so at $(x,\lambda)$ for all $x$ and for all $\lambda$ in $\epsilon$-neighborhoods of $\bar{x}$ and $\bar{\lambda}$, respectively. That is, if we do not wander too far away from $(\bar{x},\bar{\lambda})$ either in $x$ or in $\lambda$, the prox-regularity property of $f$ and its corresponding prox-regularity parameters are preserved.

\section{Examples}
\label{section:ex}
\subsection{Basic Examples}
We begin with a few basic examples of para-prox-regular functions. The first example is the result of multiplying the function $|x|$ by the parameter $\lambda$.

\begin{ex}
For $x\in\mathbb{R}^n$, $\lambda\in\mathbb{R}$, define
$$f(x,\lambda):=\lambda|x|.$$
Then $f$ is para-prox-regular at $\bar{x}=0$ with compatible parametrization by $\lambda$ at any $\bar{\lambda}>0$. It is not para-prox-regular at $\bar{x}=0$ at any $\bar{\lambda}\leq0$.
\end{ex}
\textbf{Proof:} Let $\bar{x}=0$ and $\bar{\lambda}>0$. Then\
$$f(x',\lambda)\geq f(x,\lambda)+\langle v,x'-x\rangle$$
for any $x',x$ and for any $v\in\partial_xf(x,\lambda)$, as $f$ is convex in $x$ and $\lambda>0$. Thus, for $\epsilon=\frac{\bar{\lambda}}{2}$ and $r=0$ we have
$$f(x',\lambda)\geq f(x,\lambda)+\langle v,x'-x\rangle-\frac{r}{2}|x'-x|^2$$
for $|x'-\bar{x}|<\epsilon$, $|x-\bar{x}|<\epsilon$, $|f(x,\lambda)-f(\bar{x},\bar{\lambda})|<\epsilon$, $v\in\partial_xf(x,\lambda)$, $|v-\bar{v}|<\epsilon$, and $|\lambda-\bar{\lambda}|<\epsilon$. That is, $f$ is para-prox-regular at $(\bar{x},\bar{\lambda})$.

For $\bar{\lambda}<0$, $f$ ceases to be a prox-regular function at $\bar{x}=0.$ This can be verified by noting that $\partial_xf(\bar{x},\bar{\lambda})=\{|\bar{\lambda}|,-|\bar{\lambda}|\}$ and that there does not exist $(\epsilon,r)$ that satisfies Definition \ref{df:pr1} in this case. Hence, it is not para-prox-regular either. For $\bar{\lambda}=0$ para-prox-regularity fails as well, since we require $f$ to be prox-regular in a neighborhood of $\bar{\lambda}$, but any neighborhood of $\bar{\lambda}=0$ necessarily contains some $\lambda<0$ and yields a function that is not prox-regular.\qed\\

\par The previous example hints at the next result: if a function is convex in terms of $x$ no matter the value of $\lambda$, then it is para-prox-regular. This result will be used in Subsection \ref{sub:PA} to prove that the proximal average is para-prox-regular.
\begin{ex}
\label{ex:conv}
Let $f:\mathbb{R}^n\times\mathbb{R}^m\rightarrow\mathbb{R}\cup\{\infty\}$ be convex as a function of $x$ for all $\lambda$ such that $|\lambda-\bar{\lambda}|<\epsilon$, where $\epsilon>0$. Then $f$ is para-prox-regular with compatible parametrization by $\lambda$ at $\bar{\lambda}$.
\end{ex}
\textbf{Proof:} Since $f(x,\bar{\lambda})$ is convex in terms of $x$, for any $x\in\dom f(\cdot,\bar{\lambda})$ and any $v\in\partial_xf(x,\bar{\lambda})$ we have that
\begin{equation}
f(x',\bar{\lambda})\geq f(x,\bar{\lambda})+\langle v,x'-x\rangle
\end{equation}
 for all $x'$. For every $\lambda$ such that $|\lambda-\bar{\lambda}|<\epsilon$, since $f$ remains convex, for any $x\in\dom f(\cdot,\lambda)$ and any $v_\lambda\in\partial f_x(x,\lambda)$ we have that
\begin{equation}
\label{eq:conv2}
f(x',\lambda)\geq f(x,\lambda)+\langle v_\lambda,x'-x\rangle
\end{equation}
for all $x'$. In particular, inequality (\ref{eq:conv2}) is true for all $x'$ and $x$ such that $|x'-\bar{x}|<\epsilon$, $|x-\bar{x}|<\epsilon$ and $|f(x,\lambda)-f(\bar{x},\bar{\lambda})|<\epsilon$ for any $\bar{x}$ fixed, and for all $v_\lambda\in\partial_xf(x,\lambda)$ such that $|v_\lambda-v|<\epsilon.$ This demonstrates para-prox-regularity with parameters $r=0$ and $\epsilon$.\qed
\par The next two examples demonstrate that we can take a prox-regular function $f$ and use the parameter $\lambda$ in the argument of $f$ to create para-prox-regular functions at particular $\bar{\lambda}$ values. Example \ref{ex:arg1} examines a linear shift of the argument by $\lambda$ and example \ref{ex:arg2} examines a scalar multiple of the argument by $\lambda$.
\begin{ex}
\label{ex:arg1}
Let $\tilde{f}:\mathbb{R}^n\rightarrow\mathbb{R}\cup\{\infty\}$ be prox-regular at $\bar{x}$ for $\bar{v}\in\partial\tilde{f}(\bar{x})$, and let $\lambda\in\mathbb{R}^n$. Define
$$f(x,\lambda):=\tilde{f}(x-\lambda).$$
Then $f$ is para-prox-regular at $\bar{x}$ with compatible parametrization by $\lambda$ at $\bar{\lambda}=0$.
\end{ex}
\textbf{Proof:} By assumption, there exist $\tilde{\epsilon}>0$ and $\tilde{r}\geq0$ such that
\begin{equation}
\label{eq:arg1}
\tilde{f}(\tilde{x}')\geq\tilde{f}(\tilde{x})+\langle v,\tilde{x}'-\tilde{x}\rangle-\frac{\tilde{r}}{2}|\tilde{x}'-\tilde{x}|^2
\end{equation}
whenever $\tilde{x}'\neq\tilde{x}$, $|\tilde{x}'-\bar{x}|<\tilde{\epsilon}$, $|\tilde{x}-\bar{x}|<\tilde{\epsilon}$, $|\tilde{f}(\tilde{x})-\tilde{f}(\bar{x})|<\tilde{\epsilon}$, $v\in\partial\tilde{f}(\tilde{x})$, $|v-\bar{v}|<\tilde{\epsilon}$. We need to show that when $\bar{\lambda}=0$, for all $\bar{v}\in\partial\tilde{f}(\bar{x})=\partial_xf(\bar{x},\bar{\lambda})$ there exist $\epsilon>0$ and $r\geq0$ such that
\begin{equation}
\label{eq:arg2}
f(x',\lambda)\geq f(x,\lambda)+\langle v,x'-x\rangle-\frac{r}{2}|x'-x|^2,
\end{equation}
when $x'\neq x$, $|x'-\bar{x}|<\epsilon$, $|x-\bar{x}|<\epsilon$, $|\lambda-\bar{\lambda}|<\epsilon$, $|f(x,\lambda)-f(\bar{x},\bar{\lambda})|<\epsilon$, $v\in\partial_xf(x,\lambda)$ and $|v-\bar{v}|<\epsilon$. Let $\epsilon=\frac{\tilde{\epsilon}}{2}$, $r=\tilde{r}$. Since $\epsilon<\tilde{\epsilon}$, inequality (\ref{eq:arg1}) holds when we replace $\tilde{r}$ with $r$ and $\tilde{\epsilon}$ with $\epsilon$ everywhere. Select any $\lambda$ such that $|\lambda-\bar{\lambda}|<\epsilon$, and set $\tilde{x}'=x'-\lambda$ and $\tilde{x}=x-\lambda$. Notice that if $|x'-\bar{x}|<\epsilon$ and $|x-\bar{x}|<\epsilon$, then by the triangle inequality $|(x'-\lambda)-\bar{x}|<\tilde{\epsilon}$ and $|(x-\lambda)-\bar{x}|<\tilde{\epsilon}$, which is to say $|\tilde{x}'-\bar{x}|<\tilde{\epsilon}$ and $|\tilde{x}-\bar{x}|<\tilde{\epsilon}$. If in addition we restrict our function values using $\epsilon$ instead of $\tilde{\epsilon}$, then we have
$$|\tilde{f}(\tilde{x})-\tilde{f}(\bar{x})|<\epsilon\Leftrightarrow|\tilde{f}(x-\lambda)-\tilde{f}(\bar{x}-\bar{\lambda})|<\epsilon\Leftrightarrow|f(x,\lambda)-f(\bar{x},\bar{\lambda})|<\epsilon.$$
Finally we observe that $v\in\partial\tilde{f}(\tilde{x})$ if and only if $v\in\partial_xf(x,\lambda)$, so we may rewrite inequality (\ref{eq:arg1}) as
$$\tilde{f}(x'-\lambda)\geq\tilde{f}(x-\lambda)+\langle v,(x'-\lambda)-(x-\lambda)\rangle-\frac{r}{2}|(x'-\lambda)-(x-\lambda)|^2$$
when $x'\neq x$, $|x'-\bar{x}|<\epsilon$, $|x-\bar{x}|<\epsilon$, $|\lambda-\bar{\lambda}|<\epsilon$, $|f(x,\lambda)-f(\bar{x},\bar{\lambda})|<\epsilon$, $v\in\partial_xf(x,\lambda)$, $|v-\bar{v}|<\epsilon$. That is, inequality (\ref{eq:arg2}) holds under the appropriate conditions.\qed
\begin{ex}
\label{ex:arg2}
Let $\tilde{f}:\mathbb{R}^n\rightarrow\mathbb{R}\cup\{\infty\}$ be prox-regular at $\bar{x}$ for $\bar{v}$, and let $\lambda\in\mathbb{R}$. Define
$$f(x,\lambda):=\tilde{f}(\lambda x).$$
Then $f$ is para-prox-regular at $\bar{x}$ for $\bar{v}$ with compatible parametrization by $\lambda$ at $\bar{\lambda}=1$.
\end{ex}
\textbf{Proof:} By assumption, there exist prox-regularity parameters $\tilde{\epsilon}>0$ and $\tilde{r}\geq0$ such that
\begin{equation}
\label{eq:scalar1}
\tilde{f}(\tilde{x}')\geq \tilde{f}(\tilde{x})+\langle v,\tilde{x}'-\tilde{x}\rangle+\frac{\tilde{r}}{2}|\tilde{x}'-\tilde{x}|^2
\end{equation}
for $|\tilde{x}'-\bar{x}|<\tilde{\epsilon}$, $|\tilde{x}-\bar{x}|<\tilde{\epsilon}$, $|\tilde{f}(\tilde{x})-\tilde{f}(\bar{x})|<\tilde{\epsilon}$, $v\in\partial\tilde{f}(\tilde{x})$, and $|v-\bar{v}|<\tilde{\epsilon}$. Substituting $\tilde{x}'=\lambda x'$ and $\tilde{x}=\lambda x$, inequality (\ref{eq:scalar1}) becomes
\begin{equation}
\label{eq:scalar2}
\tilde{f}(\lambda x')\geq\tilde{f}(\lambda x)+\langle v,\lambda x'-\lambda x\rangle-\frac{\tilde{r}}{2}|\lambda x'-\lambda x|^2
\end{equation}
for $|\lambda x'-\bar{x}|<\tilde{\epsilon}$, $|\lambda x-\bar{x}|<\tilde{\epsilon}$, $|\tilde{f}(\lambda x)-\tilde{f}(\bar{x})|<\tilde{\epsilon}$, $v\in\partial\tilde{f}(\lambda x)$, and $|v-\bar{v}|<\tilde{\epsilon}$. The subgradient chain rule \cite[Ex. 10.7]{rockwets} gives us $\lambda v\in\partial_xf(x,\lambda)$, so this may be rewritten as
\begin{equation}
\label{eq:scalar3}
f(x',\lambda)\geq f(x,\lambda)+\langle\lambda v,x'-x\rangle-\frac{\tilde{r}|\lambda|^2}{2}|x'-x|^2
\end{equation}
for $|\lambda x'-\bar{x}|<\tilde{\epsilon}$, $|\lambda x-\bar{x}|<\tilde{\epsilon}$, $|f(x,\lambda)-f(\bar{x},\bar{\lambda})|<\tilde{\epsilon}$, $\lambda v\in\partial_xf(x,\lambda)$, and $|v-\bar{v}|<\tilde{\epsilon}$. If we restrict our neighborhood of $\bar{x}$ further, say $|x'-\bar{x}|<\frac{\tilde{\epsilon}}{2}$ and $|x-\bar{x}|<\frac{\tilde{\epsilon}}{2}$, then there exists $\hat{\epsilon}>0$ such that for all $\lambda$ with $|\lambda-\bar{\lambda}|<\hat{\epsilon}$ and if $|x'-\bar{x}|<\frac{\tilde{\epsilon}}{2}$ and $|x-\bar{x}|<\frac{\tilde{\epsilon}}{2}$, then we necessarily have that $|\lambda x'-\bar{x}|<\tilde{\epsilon}$ and $|\lambda x-\bar{x}|<\tilde{\epsilon}$. Hence by taking $\epsilon=\min\{\frac{\tilde{\epsilon}}{2},\hat{\epsilon}\}$, inequality (\ref{eq:scalar3}) remains true for $|x'-\bar{x}|<\epsilon$, $|x-\bar{x}|<\epsilon$, $|f(x,\lambda)-f(\bar{x},\bar{\lambda})|<\epsilon$, $\lambda v\in\partial_xf(x,\lambda)$, $|v-\bar{v}|<\epsilon$, and $|\lambda-\bar{\lambda}|<\epsilon$. Since $|\lambda - \bar{\lambda}|<\epsilon$, we may take $r=\max\{\tilde{r},\tilde{r}|\bar{\lambda}+\epsilon|^2\}$. Renaming $\lambda v$ as $\tilde{v}$ we have that
\begin{equation}
f(x',\lambda)\geq f(x,\lambda)+\langle\tilde{v},x'-x\rangle-\frac{r}{2}|x'-x|^2
\end{equation}
for $|x'-\bar{x}|<\epsilon$, $|x-\bar{x}|<\epsilon$, $|\lambda-\bar{\lambda}|<\epsilon$, $|f(x,\lambda)-f(\bar{x},\bar{\lambda})|<\epsilon$, $\tilde{v}\in\partial_xf(x,\lambda)$, and $|\tilde{v}-\bar{v}|<\epsilon$.\qed

\subsection{Proximal Average}
\label{sub:PA}
The next examples are those of the \emph{proximal average} and \emph{NC-proximal average}, for which we will need some background information. Recall that for a proper function $f$, the \emph{Fenchel conjugate} $f^*$ is defined via
$$f^*(y):=\sup\limits_y\{\langle x,y\rangle-f(y)\}.$$
For two convex functions $f_0$ and $f_1$, using the Fenchel conjugate, we define the proximal average.
\begin{df}\cite[Definition 4.1]{howto}
The \emph{proximal average} $PA_{f_0,f_1}$ of two proper convex functions $f_0$ and $f_1$ is defined via
$$PA_{f_0,f_1}(x,\lambda):=((1-\lambda)(f_0+\frac{1}{2}|\cdot|^2)^*+\lambda(f_1+\frac{1}{2}|\cdot|^2)^*)^*(x)-\frac{1}{2}|x|^2$$
for $\lambda\in[0,1]$.\end{df}
The proximal average has been shown to be an effective method of transforming one convex function into another in a continuous manner \cite{howto}. It has also been shown to be a convex function \cite{proxpoint}, and as such it is prox-regular. In fact, it is para-prox-regular, as the following proposition demonstrates.
\begin{prop}
The proximal average $PA_{f_0,f_1}$ of two proper convex functions $f_0,f_1$ is para-prox-regular with compatible parametrization by $\lambda$ at $\bar{\lambda}$, for all $\bar{\lambda}\in(0,1)$.
\end{prop}
\textbf{Proof:} Bauschke, Matou\v{s}kov\'a and Reich \cite[Theorem 6.1]{proxpoint} proved that this function is convex in $x$ for all $\lambda\in(0,1)$, therefore it satisfies the conditions of Example \ref{ex:conv}.\qed\medskip\\
In \cite{proxave} an alternate representation of the proximal average function was given using Moreau envelopes. Recall that for a function $f$ the \emph{Moreau envelope} with parameter $r>0$ is defined as
$$e_rf(x):=\inf\limits_y\{f(y)+\frac{r}{2}|y-x|^2\},$$
and its associated \emph{proximal point mapping} is defined as
$$P_rf(x):=\argmin\limits_y\{f(y)+\frac{r}{2}|y-x|^2\}.$$
If $f$ is proper and convex, then $e_rf$ is proper for all $r>0$. Properties of the Moreau envelope give rise to the equivalent definition of proximal average \cite[Equation (4)]{proxave},
$$PA_{f_0,f_1}(x,\lambda)=-e_1(-(1-\lambda)e_1f_0-\lambda e_1f_1)(x).$$\\
However, this formulation is only useful if $f_0$ and $f_1$ are convex functions. For two not necessarily convex functions $f_0$ and $f_1$, the \emph{NC-proximal average function} $PA_r$ is similarly defined. We need one more definition first, that of prox-boundedness.
\begin{df}
A function $f$ is \emph{prox-bounded} if there exist $r>0$ and a point $\bar{x}$ such that $e_rf(\bar{x})>-\infty$. The infimum of the set of all such $r$ is called the \emph{threshold of prox-boundedness}.
\end{df}
\begin{df}\cite[Equation (5)]{proxave}
Let $f_0$ and $f_1$ be proper, prox-bounded functions, and let $r>0$ be greater than the threshold of prox-boundedness for both $f_0$ and $f_1$. The \emph{NC-proximal average} of $f_0$ and $f_1$ is defined as
$$PA_r(x,\lambda):=-e_{r+\lambda(1-\lambda)}(-(1-\lambda)e_rf_0-\lambda e_rf_1)(x).$$
\end{df}

Showing that the NC-proximal average is para-prox-regular requires a technical condition from \cite[Theorem 4.6]{proxave}, specifically that the Lipschitz constant of $\lambda P_rf_0+(1-\lambda)P_rf_1-I$ is bounded by $1$.  This condition ensures that the NC-proximal average is well-behaved.

\begin{ex}Consider two proper, prox-bounded, prox-regular functions $f_0$ and $f_1$.
If the Lipschitz constant of $\lambda P_rf_0+(1-\lambda)P_rf_1-I$ is bounded by $1$, then for $r$ sufficiently large, the NC-proximal average is para-prox-regular for all $\bar{\lambda}\in(0,1)$.
\end{ex}
\textbf{Proof:} Fix $\bar{\lambda}\in(0,1)$. Let $\bar{v}\in\partial_xPA_r(\bar{x},\bar{\lambda})$. Enlarging $r$ if necessary, we know that $PA_r$ is locally Lipschitz continuous in $\lambda$ \cite[Theorem 4.6]{proxave}, and that it is a lower-$\mathcal{C}^2$ function in $x$ \cite[Proposition 2.5]{proxave}. Since all lower-$\mathcal{C}^2$ functions are prox-regular, we have that $PA_r$ is prox-regular as a function of $x$ at $\bar{x}$ for $\bar{v}$, say with parameters $\epsilon_1>0$ and $\rho\geq0$. Therefore
\begin{equation}
\label{eq:pa1}
PA_r(x',\bar{\lambda})>PA_r(x,\bar{\lambda})+\langle v,x'-x\rangle-\frac{\rho}{2}|x'-x|^2
\end{equation}
when $x'\neq x$, $|x'-\bar{x}|<\epsilon_1$, $|x-\bar{x}|<\epsilon_1$, $|PA_r(x,\bar{\lambda})-PA_r(\bar{x},\bar{\lambda})|<\epsilon_1$, $v\in\partial_xPA_r(x,\bar{\lambda})$, and $|v-\bar{v}|<\epsilon_1$. We have established that $PA_r$ is continuous in $\lambda$, and therefore there exists $\epsilon_2>0$ such that inequality (\ref{eq:pa1}) remains true when $\bar{\lambda}$ is replaced by any $\lambda$ such that $|\lambda-\bar{\lambda}|<\epsilon_2$. We can reduce $\epsilon_2$ if necessary, so that $|\lambda-\bar{\lambda}|<\epsilon_2$ implies $\lambda\in(0,1)$. Taking $\epsilon_3=\min\{\epsilon_1,\epsilon_2\}$, we have
\begin{equation}
\label{eq:pa2}
PA_r(x',\lambda)>PA_r(x,\lambda)+\langle v,x'-x\rangle-\frac{\rho}{2}|x'-x|^2
\end{equation}
when $x'\neq x$, $|x'-\bar{x}|<\epsilon_3$, $|x-\bar{x}|<\epsilon_3$, $|PA_r(x,\lambda)-PA_r(\bar{x},\bar{\lambda})|<\epsilon_3$, $v\in\partial_xPA_r(x,\bar{\lambda})$, $|v-\bar{v}|<\epsilon_3$, and $|\lambda-\bar{\lambda}|<\epsilon_3$. If we use $\frac{\epsilon_3}{2}$-neighborhoods instead of $\epsilon_3$-neighborhoods, then inequality (\ref{eq:pa2}) certainly remains true. Since $PA_r$ is $\mathcal{C}^1$ with respect to $x,$ we know that $\partial_xPA_r(x,\bar{\lambda})=\{\nabla_xPA_r(x,\bar{\lambda})\},$ and by \cite[Theorem 4.6]{proxave} we know that $\nabla_xPA_r(x,\bar{\lambda})$ changes in a continuous manner with respect to $\lambda.$ Hence, there exists $\epsilon_4>0$ such that inequality (\ref{eq:pa2}) remains valid when $v\in\partial_xPA_r(x,\bar{\lambda})$ is replaced by $\hat{v}\in\partial_xPA_r(x,\lambda)$ for $|\hat{v}-v|<\epsilon_4$ and $|\lambda-\bar{\lambda}|<\epsilon_4$. Shrinking $\epsilon_4$ if necessary, we can ensure that $\epsilon_4<\frac{\epsilon_3}{2}$. Thus we have $|\hat{v}-v|<\frac{\epsilon_3}{2}$ and $|v-\bar{v}|<\frac{\epsilon_3}{2}$, which implies that $|\hat{v}-\bar{v}|<\frac{\epsilon_3}{2}+\frac{\epsilon_3}{2}=\epsilon_3$. Therefore, choosing $\epsilon=\min\{\frac{\epsilon_3}{2},\epsilon_4\}$ we have that
\begin{equation}
\label{eq:pa3}
PA_r(x',\lambda)>PA_r(x,\lambda)+\langle \hat{v},x'-x\rangle-\frac{\rho}{2}|x'-x|^2
\end{equation}
when $x'\neq x$, $|x'-\bar{x}|<\epsilon$, $|x-\bar{x}|<\epsilon$, $|PA_r(x,\lambda)-PA_r(\bar{x},\bar{\lambda})|<\epsilon$, $\hat{v}\in\partial_xPA_r(x,\lambda)$, $|\hat{v}-\bar{v}|<\epsilon$, and $|\lambda-\bar{\lambda}|<\epsilon$.\qed

\section{Monotonicity}
\label{section:mono}
This section contains an alternate definition of para-prox-regular functions. We begin with the definition of an \emph{$f$-attentive localization} of the subdifferential. An $f$-attentive localization of $\partial_xf$ around $(\bar{x}, \bar{\lambda}, \bar{v})$ for a fixed $\lambda$ is a set-valued mapping $T_\lambda:\mathbb{R}^n\rightrightarrows\mathbb{R}^n$ whose graph in $\mathbb{R}^n\times\mathbb{R}^n$ is the intersection of $\gph\partial_xf$ with an $f$-attentive neighborhood of $\bar{x}$ and an ordinary neighborhood of $\bar{v}$. (This contrasts with an ordinary localization, in which the $f$-attentive neighborhood of $\bar{x}$ is relaxed to an ordinary neighborhood.) The following definition is adapted from \cite[Def. 3.1]{proxfunc}.
\begin{df}
\label{df:attent}
For an $\epsilon>0$ and a fixed $\lambda$, the \emph{$f$-attentive $\epsilon$-localization of $\partial f_x$ around $(\bar{x},\bar{\lambda},\bar{v})$} is
\begin{equation}
\label{equation:attent}
T_\lambda(x)=\begin{cases}
\{v\in\partial_xf(x,\lambda):|v-\bar{v}|<\epsilon\}, &  |x-\bar{x}|<\epsilon, |f(x,\lambda)-f(\bar{x},\bar{\lambda})|<\epsilon\\
\emptyset, & otherwise.
\end{cases}
\end{equation}
\end{df}
With this definition in mind, we shall present Theorem \ref{thm:mono}, which is a para-prox-regular analog of \cite[Theorem 3.2]{proxfunc}, and Corollary \ref{cor:alternate}, which is the alternate definition of para-prox-regular functions that we seek. In order to prove Theorem \ref{thm:mono}, we use the following lemma.  The lemma takes a para-prox-regular function at $\bar{x}$ for $\bar{v}$ (with compatible parametrization by $\lambda$ at $\bar{\lambda}$) and forms a new function that is para-prox-regular at $0$ for $0$ (with compatible parametrization by $\lambda$ at $\bar{\lambda}$), in effect shifting the original function. This will allow us, in Theorem \ref{thm:mono}, to make some assumptions without loss of generality in order to simplify the proof.

\begin{lem}
\label{lem:vbar}
Let $f$ be para-prox-regular at $\bar{x}$ for $\bar{v}$ with compatible parametrization by $\lambda$ at $\bar{\lambda}$. Define
$$\tilde{f}(x,\lambda):=f(x-\bar{x},\lambda)-\langle\bar{v},x-\bar{x}\rangle.$$
Then $\tilde{f}$ is para-prox-regular at $0$ for $0$ with compatible parametrization by $\lambda$ at $\bar{\lambda}$.
\end{lem}
\textbf{Proof:} By assumption, there exist $\epsilon>0$ and $r\geq0$ such that
\begin{equation}
\label{eq:paralem}
f(x',\lambda)\geq f(x,\lambda)+\langle v,x'-x\rangle-\frac{r}{2}|x-\bar{x}|^2
\end{equation}
whenever $x'\neq x$, $|x'-\bar{x}|<\epsilon$, $|x-\bar{x}|<\epsilon$, $|f(x,\lambda)-f(\bar{x},\bar{\lambda})|<\epsilon$, $v\in\partial_xf(x,\lambda)$, $|v-\bar{v}|<\epsilon$, and$|\lambda-\bar{\lambda}|<\epsilon$. For ease of discussion let $f_1(x,\lambda)=f(x-\bar{x},\lambda)$ and $f_2(x,\lambda)=-\langle\bar{v},x-\bar{x}\rangle$. Notice that $v\in\partial_xf(x,\lambda)$ implies $v-\bar{v}\in\partial_xf(x,\lambda)-\bar{v}$. Since $f_2 \in \mathcal{C}^2$, applying \cite[Exercise 8.8(c)]{rockwets} we have that
    \[\partial_xf(x,\lambda)-\bar{v}=
    \partial_x[f(x-\bar{x},\lambda)-\langle\bar{v},x-\bar{x}\rangle]=
    \partial_x\tilde{f}(x,\lambda).\]
Therefore $v-\bar{v}\in\partial_x\tilde{f}(x,\lambda)$. Translating inequality (\ref{eq:paralem}) by $\bar{x}$, we may rewrite it as
\begin{equation}
\label{equation:paralem2}
f(x'-\bar{x},\lambda)\geq f(x-\bar{x},\lambda)+\langle v,(x'-\bar{x})-(x-\bar{x})\rangle-\frac{r}{2}|(x'-\bar{x})-(x-\bar{x})|^2
\end{equation}
now with the conditions $x'\neq x$, $|x'|<\epsilon$, $|x|<\epsilon$, $|\tilde{f}(x,\lambda)-\tilde{f}(0,\bar{\lambda})|<\epsilon$, $v-\bar{v}\in\partial_x\tilde{f}(x,\lambda)$, $|v-\bar{v}|<\epsilon$, and $|\lambda-\bar{\lambda}|<\epsilon$. Subtracting $\langle\bar{v},x-\bar{x}\rangle$ from both sides of inequality (\ref{equation:paralem2}) and simplifying, we get
\begin{equation}
\label{equation:paralem3}
\tilde{f}(x',\lambda)\geq\tilde{f}(x,\lambda)+\langle v,x'-x\rangle-\frac{r}{2}|x'-x|^2.
\end{equation}
Therefore $\tilde{f}$ is para-prox-regular at $0$ for $0$ with compatible parametrization by $\lambda$ at $\bar{\lambda}$.\qed

\begin{thm}
\label{thm:mono}
Let $f:\mathbb{R}^n\times\mathbb{R}^m\rightarrow\mathbb{R}\cup\{\infty\}$ be a proper function that is locally lsc near $(\bar{x},\bar{\lambda})$. Then the following are equivalent.
\begin{itemize}
\item[a)] The function $f$ is para-prox-regular at $\bar{x}$ for $\bar{v}$ with compatible parametrization by $\lambda$ at $\bar{\lambda}$.
\item[b)] The vector $\bar{v}$ is a proximal subgradient of $f$ with respect to $x$ at $(\bar{x},\bar{\lambda})$, and there exist $\epsilon>0$, $r>0$ such that the $f$-attentive $\epsilon$-localization $T_\lambda$ of $\partial_xf$ at $(\bar{x},\bar{\lambda},\bar{v})$ has $T_\lambda+rI_x$ monotone, for all $\lambda$ such that $|\lambda-\bar{\lambda}|<\epsilon$. That is, there is a constant $r$ such that
\begin{equation}
\label{equation:mono}
\langle v_1-v_0,x_1-x_0\rangle\geq-r|x_1-x_0|^2
\end{equation}
when $v_i\in T_\lambda(x_i)$, $i\in\{0,1\}$, and $|\lambda-\bar{\lambda}|<\epsilon$. Here $I_x$ is the identity mapping with respect to $x$.
\item[c)] There exist $\epsilon>0, r>0$ such that
\begin{equation}
\label{equation:mono4}
f(x,\lambda)>f(\bar{x},\bar{\lambda})+\langle\bar{v},x-\bar{x}\rangle-\frac{r}{2}|x-\bar{x}|^2
\end{equation}
when $|x-\bar{x}|<\epsilon$ and $|\lambda-\bar{\lambda}|<\epsilon$, and the mapping $(\partial_xf+rI_x)^{-1}$ has the following single-valuedness property near $\bar{z}=\bar{v}+r\bar{x}$: if $|z-\bar{z}|<\epsilon$, then for $i\in\{0,1\}$ one has
\begin{equation}
\label{equation:mono5}
(x_i,\lambda)\in(\partial_xf+rI_x)^{-1}(z)\Rightarrow x_0=x_1
\end{equation}
when $|x_i-\bar{x}|<\epsilon$, $|f(x_i,\lambda)-f(\bar{x},\bar{\lambda})|<\epsilon$, and $|\lambda-\bar{\lambda}|<\epsilon$.
\end{itemize}
\end{thm}
\textbf{Proof:}\\
\smallskip
a)$\Rightarrow$b): We have that $f$ is para-prox-regular at $(\bar{x},\bar{\lambda})$ for $\bar{v}$, so there exist $\epsilon>0$ and $r>0$ such that
\begin{equation}
\label{equation:mono3}
\begin{array}{c}
f(x_1,\lambda)\geq f(x_0,\lambda)+\langle v_0,x_1-x_0\rangle-\frac{r}{2}|x_1-x_0|^2\\
f(x_0,\lambda)\geq f(x_1,\lambda)+\langle v_1,x_0-x_1\rangle-\frac{r}{2}|x_0-x_1|^2
\end{array}
\end{equation}
whenever $i\in\{0,1\}$, $|x_i-\bar{x}|<\epsilon$, $|f(x_i,\lambda)-f(\bar{x},\bar{\lambda})|<\epsilon$, $v_i\in\partial_xf(x_i,\lambda)$, $|v_i-\bar{v}|<\epsilon$, and $|\lambda-\bar{\lambda}|<\epsilon.$ Adding the two inequalities in (\ref{equation:mono3}) together, we get
$$f(x_1,\lambda)+f(x_0,\lambda)\geq f(x_0,\lambda)+f(x_1,\lambda)+\langle v_0-v_1,x_1-x_0\rangle-r|x_1-x_0|^2,$$
which simplifies to inequality (\ref{equation:mono}). Finally, since $f$ is para-prox-regular at $\bar{x}$ for $\bar{v}$, it is prox-regular there as well, giving us that $\bar{v}$ is a proximal subgradient.\medskip\\
b)$\Rightarrow$c): For this section of the proof and the next, without loss of generality we assume that $\bar{x}=0$, $f(0,\bar{\lambda})=0$, and $\bar{v}=0$ (by Lemma \ref{lem:vbar}). Let $\bar{\epsilon}>0$ and $\bar{r}>0$ be the parameters assumed by condition b). Since $\bar{v}$ is a proximal subgradient of $f$ with respect to $x$ at $(\bar{x},\bar{\lambda})$, we have
\begin{equation}
\label{equation:mono6}
f(x,\lambda)>f(\bar{x},\bar{\lambda})+\langle\bar{v},x-\bar{x}\rangle-\frac{r}{2}|x-\bar{x}|^2
\end{equation}
when $|x-\bar{x}|<\epsilon$ and $|\lambda-\bar{\lambda}|<\epsilon$, for some $\epsilon\in(0,\bar{\epsilon})$ and some $r\in(\bar{r},\infty)$. Let $\lambda\in[\bar{\lambda}-\epsilon,\bar{\lambda}+\epsilon]$. By decreasing $\epsilon$ and increasing $r$ if necessary, we can arrange that $r\geq1$ and $\epsilon<\frac{\bar{\epsilon}}{2}$. Suppose that for $z$ such that $|z-(\bar{v}+r\bar{x})|<\epsilon$ we have $(x_i,\lambda)\in(\partial_xf+rI_x)^{-1}(z)$, $i\in\{0,1\}$. By our assumptions at the beginning of the proof, this simplifies to $|z|<\epsilon$. For $v_i=z-rx_i$ we have $v_i\in\partial_xf(x_i,\lambda)$, and by the triangle inequality $|v_i|\leq|z|+r|x_i|<\frac{\bar{\epsilon}}{2}+\frac{\bar{\epsilon}}{2}=\bar{\epsilon}$. Also $|x_i|<\bar{\epsilon}$ and $|f(x_i,\lambda)|<\bar{\epsilon}$, thus, $v_i\in T_\lambda(x_i)$. So we have that
$$-\bar{r}|x_1-x_0|^2\leq\langle v_1-v_0,x_1-x_0\rangle=\langle(z-rx_1)-(z-rx_0),x_1-x_0\rangle=-r|x_1-x_0|^2.$$
Since $r>\bar{r}$, we conclude that $x_1=x_0$.\medskip\\
c)$\Rightarrow$a): Let $\epsilon>0$ and $r>0$ be the parameters assumed by condition c). It will suffice to show that there exists $\bar{\epsilon}\in(0,\epsilon)$ such that para-prox-regularity of $f$ is satisfied by
$\bar{\epsilon}$ and $r$ at $\bar{x}=0$, $f(0,\bar{\lambda})=0$ for $\bar{v}=0$. That is, whenever
\begin{equation}
\label{equation:mono8}
v_0\in\partial_xf(x_0,\lambda), |v_0|<\bar{\epsilon}, |f(x_0,\lambda)|<\bar{\epsilon},
\end{equation}
we have
$$f(x,\lambda)\geq f(x_0,\lambda)+\langle v_0,x-x_0\rangle-\frac{r}{2}|x-x_0|^2$$
when $|x|<\bar{\epsilon}$ and $|\lambda-\bar{\lambda}|<\bar{\epsilon}$. In place of the latter we can aim at guaranteeing the stronger condition
\begin{equation}
\label{equation:mono9}
f(x,\lambda)>f(x_0,\lambda)+\langle v_0,x-x_0\rangle-\frac{r}{2}|x-x_0|^2
\end{equation}
for $x\neq x_0$, $|x|\leq\epsilon$, and $|\lambda-\bar{\lambda}|\leq\epsilon$.
Notice that
$$-\frac{r}{2}|x-x_0|^2=\langle rx_0,x-x_0\rangle-\frac{r}{2}|x|^2+\frac{r}{2}|x_0|^2,$$
so we can rewrite inequality (\ref{equation:mono9}) as
\begin{align*}
 f(x,\lambda) & >f(x_0,\lambda)+\langle v_0,x-x_0\rangle+\langle rx_0,x-x_0\rangle-\frac{r}{2}|x|^2+\frac{r}{2}|x_0|^2\\
& =f(x_0,\lambda)+\langle v_0+rx_0,x\rangle-\langle v_0+rx_0,x_0\rangle-\frac{r}{2}|x|^2+\frac{r}{2}|x_0|^2.
\end{align*}
Rearranging, we get
$$f(x,\lambda)-\langle v_0+rx_0,x\rangle+\frac{r}{2}|x|^2>f(x_0,\lambda)-\langle v_0+rx_0,x_0\rangle+\frac{r}{2}|x_0|^2.$$
Thus, we seek to show that for every fixed $\lambda$ such that $|\lambda-\bar{\lambda}|\leq\bar{\epsilon}$,
\begin{equation}
\label{equation:mono10}
\argmin\limits_{|x|\leq\bar{\epsilon}}\{f(x,\lambda)+\frac{r}{2}|x|^2-\langle z_0,x\rangle\}=\{x_0\},
\end{equation}
where $z_0=v_0+rx_0$. In summary, we need only demonstrate that conditions (\ref{equation:mono8}) imply equation (\ref{equation:mono10}) when $\bar{\epsilon}$ is small enough. Our assumption that $f$ is locally lsc at $(0,0)$ with $f(0,0)=0$ ensures that $f$ is lsc with respect to a compact set of the form
$$C=\{(x,\lambda):|x|\leq\hat{\epsilon}, |\lambda|\leq\hat{\epsilon}, f(x,\lambda)\leq2\hat{\epsilon}\}$$
for some $\hat{\epsilon}>0$. In particular, $f$ must be lsc at $(0,0)$ itself: $\liminf\limits_{x,\lambda\rightarrow0}f(x,\lambda)=0$.\\
\smallskip
Shrinking $\epsilon$  if necessary, we can arrange that
\begin{equation}
\label{equation:mono11}
|x|\leq\epsilon, |\lambda|\leq\epsilon\Rightarrow|x|\leq\hat{\epsilon},|\lambda|\leq\hat{\epsilon},\mbox{ and }f(x,\lambda)>-2\hat{\epsilon}.
\end{equation}
Enlarging $r$ if necessary, we can arrange that
\begin{equation}
\label{equation:mono12}
|x|^2<\frac{6\hat{\epsilon}}{r}\Rightarrow|x|<\epsilon,|x|<\frac{\epsilon^2}{2\hat{\epsilon}},\mbox{ and }f(x,\lambda)>-\frac{\epsilon}{2}.
\end{equation}
Define
$$g(z):=\inf\limits_{|x|,|\lambda-\bar{\lambda}|\leq\epsilon}\{f(x,\lambda)+\frac{r}{2}|x|^2-\langle z,x\rangle\},$$
$$G(z):=\argmin\limits_{|x|,|\lambda-\bar{\lambda}|\leq\epsilon}\{f(x,\lambda)+\frac{r}{2}|x|^2-\langle z,x\rangle\}.$$
Since $g$ is the pointwise infimum of a collection of affine functions of $z$, it is concave. (Note: $\bar{x}=0$, $\bar{v}=0$, and $f(\bar{x},\bar{\lambda})=0$.) Our assumptions tell us that $g(0)=0$ and $G(0)=\{0\}$, whereas $g(z)\leq0$ in general. We claim that under these circumstances,
\begin{equation}
\label{equation:mono13}
|z|<\frac{\hat{\epsilon}}{\epsilon}\Rightarrow\begin{cases}
G(z)\neq\emptyset\mbox{ and} & \\|f(x,\lambda)|<2\hat{\epsilon} &\mbox{ for all }(x,\lambda)\in G(z).\end{cases}
\end{equation}
To see this, observe that because $g\leq0$, the minimization in the definition of $g(z)$ is unaffected if the attention is restricted to the points $(x,\lambda)$ satisfying not only $|x|<\epsilon$ and $|\lambda|<\epsilon$, but also $f(x,\lambda)+\frac{r}{2}|x|^2-\langle z,x\rangle\leq\hat{\epsilon}$, in which case $f(x,\lambda)\leq\hat{\epsilon}+\epsilon|z|$. Therefore, as long as $|z|<\frac{\hat{\epsilon}}{\epsilon}$, attention can be restricted to points $(x,\lambda)$ satisfying $f(x,\lambda)<2\hat{\epsilon}$. Recalling conditions (\ref{equation:mono11}) and the choice of the set $C$, we deduce that when $|z|\leq\frac{\hat{\epsilon}}{\epsilon}$,
\begin{equation}
\label{equation:mono14}
\begin{array}{c c}
g(z)=\inf\limits_{(x,\lambda)\in C}\{f(x,\lambda)+\frac{r}{2}|x|^2-\langle z,x\rangle\}\mbox{ and}\\
G(z)=\argmin\limits_{(x,\lambda)\in C}\{f(x,\lambda)+\frac{r}{2}|x|^2-\langle z,x\rangle\}.
\end{array}
\end{equation}
Since $f$ is lsc relative to $C$, which is compact, the infimum is sure to be attained. Thus, implication (\ref{equation:mono13}) is correct. Also we observe that $g$ is finite on a neighborhood of $z=0$, and hence continuous (by concavity) \cite[Theorem 2.35]{rockwets}. Choose $\tilde{\epsilon}\in(0,\epsilon)$ small enough that $\tilde{\epsilon}(1+r)<\frac{\hat{\epsilon}}{\epsilon}$ and $g(z)>-\frac{\epsilon}{2}$ when $|z|<\tilde{\epsilon}(1+r)$. Under this choice we are ready to consider elements satisfying conditions (\ref{equation:mono8}) and show that we get equation (\ref{equation:mono10}). The vector $z_0=v_0+rx_0$ has $|z_0|\leq|v_0|+r|x_0|<\tilde{\epsilon}(1+r)$, hence $|z_0|<\frac{\hat{\epsilon}}{\epsilon}$ and $g(z_0)>-\frac{\epsilon}{2}$. In particular, $G(z_0)$ must be nonempty. Consider any $(x_1,\lambda)\in G(z_0)$. We have to show on the basis of the single-valuedness property that $x_1=x_0$. We have that $(x_0,\lambda)\in(\partial_xf+rI)^{-1}(z_0)$, due to the fact that $v_0\in\partial_xf(x_0,\lambda)$ and hence $v_0+rx_0\in(\partial_xf+rI_x)(x_0,\lambda)$. If we can establish that
\begin{equation}
\label{equation:mono15}
|x_1|<\epsilon, |\lambda|<\epsilon, |f(x_1,\lambda)|<\epsilon\mbox{ and }(x_1,\lambda)\in(\partial_f+rI)^{-1}(z_0),
\end{equation}
then the single-valuedness property can be invoked and we will get $x_1=x_0$ as required.\\
\smallskip
 Because $|z_0|<\frac{\hat{\epsilon}}{\epsilon}$, we know from implication (\ref{equation:mono13}) that $|f(x_1,\lambda)|<2\hat{\epsilon}$. Then from $\frac{r}{2}|x_1|^2=\langle z_0,x_1\rangle-f(x_1,\lambda)+g(z_0)$, where $g(z_0)\leq0$, we know that $\frac{r}{2}|x_1|^2\leq|z_0||x_1|+|f(x_1,\lambda)|<\frac{\hat{\epsilon}}{\epsilon}\epsilon+2\hat{\epsilon}=3\hat{\epsilon}$, so $|x_1|^2<\frac{6\hat{\epsilon}}{r}$. Through (\ref{equation:mono12}) this implies that $|x_1|<\epsilon$ and also that $|x_1|<\frac{\epsilon^2}{2\hat{\epsilon}}$ and $f(x_1,\lambda)>-\epsilon$. Then from $f(x_1,\lambda)=\langle z_0,x_1\rangle-\frac{r}{2}|x_1|^2+g(z_0)\leq|z_0||x_1|+|g(z_0)|$ and the fact that $g(z_0)<\frac{\epsilon}{2}$ we obtain $f(x_1)<\frac{\hat{\epsilon}}{\epsilon}\frac{\epsilon^2}{2\hat{\epsilon}}+\frac{\epsilon}{2}=\epsilon$. Hence $|x_1|<\epsilon$, $|\lambda|<\epsilon$, and $|f(x_1,\lambda)|<\epsilon.$ The fact that the minimum for $g(z_0)$ is attained at $x_1$ gives us
\begin{equation}
\label{equation:mono16}
f(x,\lambda)+\frac{r}{2}|x|^2-\langle z_0,x\rangle\geq f(x_1,\lambda)+\frac{r}{2}|x_1|^2-\langle z_0,x_1\rangle
\end{equation}
when $|x|<\epsilon$ and $|\lambda|<\epsilon$. Hence, for all $|x|<\epsilon$ and $|\lambda|<\epsilon,$
\begin{align*}
f(x,\lambda) & \geq f(x_1,\lambda)+\langle z_0,x\rangle-\langle z_0,x_1\rangle-\frac{r}{2}|x|^2+\frac{r}{2}|x_1|^2\\
 & =f(x_1,\lambda)+\langle z_0,x-x_1\rangle-r\langle x,x_1\rangle+r\langle x,x_1\rangle+r\langle x_1,x_1\rangle-\frac{r}{2}\langle x,x\rangle-\frac{r}{2}\langle x_1,x_1\rangle\\
 & =f(x_1,\lambda)+\langle z_0,x-x_1\rangle-\langle rx_1,x-x_1\rangle-\frac{r}{2}(\langle x,x\rangle-2\langle x,x_1\rangle+\langle x_1,x_1\rangle)\\
 & =f(x_1,\lambda)+\langle z_0-rx_1,x-x_1\rangle-\frac{r}{2}|x-x_1|^2.\\
\end{align*}
Since $|x_1|<\epsilon$, this implies that $z_0-rx_1\in\partial_x f(x_1,\lambda)$, so $z_0\in(\partial_xf+rI_x)(x_1,\lambda)$. The subgradient relation in (\ref{equation:mono15}) is therefore correct as well, and our single-valuedness assumption implies $x_0=x_1$. Thus equation (\ref{equation:mono10}) holds and we are done.\qed\\

This result allows us to define the following alternate definition of para-prox-regularity.
\begin{cor}
\label{cor:alternate}
Let $f$ be a proper, lsc function. Then $f$ is para-prox-regular in $x$ at $\bar{x}$ for $\bar{v}$ with compatible parametrization by $\lambda$ at $\bar{\lambda}$, with parameters $\epsilon>0, r>0$, if and only if $\bar{v}$ is a proximal subgradient of $f$ with respect to $x$ at $(\bar{x},\bar{\lambda})$ and
$$\langle v_1-v_0,x_1-x_0\rangle\geq r|x_1-x_0|^2$$
when for $i\in\{0,1\}$ we have $|x_i-\bar{x}|<\epsilon$, $|f(x_i,\lambda)-f(\bar{x},\bar{\lambda})|<\epsilon$, $v_i\in\partial_xf(x_i,\lambda)$, $|v_i-\bar{v}|<\epsilon$, and $|\lambda-\bar{\lambda}|<\epsilon$.
\end{cor}

\section{Advanced Examples}
\label{section:main}
In this section we examine some common situations where parameters arise in optimization problems. In all cases we focus on identifying the prox-regularity parameters. We already know of many classes of functions that are prox-regular, some were mentioned in the introduction, and here in general we begin with prox-regular functions and build upon them. We show that the scalar multiplication of a prox-regular function by a positive parameter results in a function that remains prox-regular. This will lead to the construction of para-prox-regular functions via the weighted sums and weighted max of prox-regular functions.
\par The following proposition and theorem is a result first provided in \cite{calcprox}, showing that strongly amenable functions are prox-regular. In that work, Poliquin and Rockafellar simply conclude the function is prox-regular and they do not single out the prox-regularity parameters $\epsilon$ and $r$ in the statement of the theorem. In the present work, the theorem has been restated with a proposition beforehand, in order to pull out the prox-regularity parameters from the proof into the statement of the theorem. In later results we require the ability to identify such parameters instead of merely stating that they exist. Other than clarifying the prox-regularity parameters, Proposition \ref{prop:strong} and Theorem \ref{thm:prelim} and their proofs are essentially the same as \cite[Theorem 3.1]{calcprox}.
\begin{prop}
\label{prop:strong}
Let $f(x)=g(F(x))$, where $F:\mathbb{R}^n\rightarrow\mathbb{R}^m$ is of class $\mathcal{C}^2, g:\mathbb{R}^m\rightarrow\bar{\mathbb{R}}$ is lsc and proper, and suppose there is no vector $y\neq0$ in $\partial^\infty g(F(\bar{x}))$ with $\nabla F(\bar{x})^*y=0$, where $\bar{x}\in\dom f$. Then
\begin{itemize}
\item[a)] there exists $\epsilon_0>0$ such that
$$\partial f(x)\subset\nabla F(x)^*\partial g(F(x))$$
when $|x-\bar{x}|<\epsilon_0$ and $|f(x)-f(\bar{x})|<\epsilon_0$.
\item[b)]$S:(x,v)\rightarrow\{y\in\partial g(F(x)):\nabla F(x)^*y=v\}$
has closed graph and is locally bounded at $(\bar{x},\bar{v})$ for $\bar{v}\in\partial f(\bar{x})$, within the same neighborhoods $B_{\epsilon_0}(\bar{x})$ and $B_{\epsilon_0}(f(\bar{x}))$. In particular, $S(\bar{x},\bar{v})$ is a compact set.
\item[c)] If $\bar{v}\in\partial f(\bar{x})$ is such that for every $y\in\partial g(F(\bar{x}))$ with $\nabla F(\bar{x})^*y=\bar{v}$ the function $g$ is prox-regular at $F(\bar{x})$ with parameters $\epsilon_y>0, r_y>0$, then for the parameters $\epsilon=\min\{\epsilon_0,\min\limits_y\{\epsilon_y\}\}$ and $\bar{r}=\max\limits_y\{r_y\}$, which are guaranteed to exist since $S(\bar{x},\bar{v})$ is a compact set, we have
\begin{equation}
\label{eq:stat2}
\langle y_1-y_0,u_1-u_0\rangle\geq-\bar{r}|u_1-u_0|^2
\end{equation}
whenever $y_i\in\partial g(u_i), |u_i-F(\bar{x})|<\bar{\epsilon}, |g(u_i)-g(F(\bar{x}))|<\bar{\epsilon}, \dist(y_i,S(\bar{x},\bar{v}))<\bar{\epsilon}$.\\
\end{itemize}
\end{prop}
\textbf{Proof:} See \cite[Theorem 3.1]{calcprox} and the proof thereof.\medskip\\
In Theorem \ref{thm:prelim} we must provide greater detail, in order to fully describe the prox-regularity parameters. In particular, the full details of parameter $r$ only become apparent with the proof itself.

\begin{thm}
\label{thm:prelim}
Let the conditions of Proposition \ref{prop:strong} hold. Assume further that $\bar{v}\in\partial f(\bar{x})$ is a vector such that for every $y\in\partial g(F(\bar{x}))$ such that $\nabla F(\bar{x})^*y=\bar{v}$, the function $g$ is prox-regular at $F(\bar{x})$ for $y$ with parameters $\epsilon_y>0, r_y>0$. Then $f$ is prox-regular at $\bar{x}$ for $\bar{v}$ with parameters $\epsilon>0$ as in Proposition \ref{prop:strong} and $r>0$ as described by equations (\ref{eqrbar}), (\ref{eqr1}), (\ref{eqr2}), and (\ref{eqr}) below.
\end{thm}
\medskip
\textbf{Proof:} By our assumptions, $g$ is prox-regular at $F(\bar{x})$ for each $y\in S(\bar{x},\bar{v})$, with constants $\epsilon_y$ and $r_y$. Since $S(\bar{x},\bar{v})$ is compact we can make this uniform. We know that there is a finite open cover of $S(\bar{x},\bar{v})$ using a finite number of the $\epsilon_y,$ say $\{\epsilon_1,\ldots,\epsilon_k\}.$ Then using the corresponding $\{r_1,\ldots,r_k\},$ we define
	\begin{equation}\bar{r}:=\max_{i\in\{1,\ldots,k\}}\{r_i\}>0.\label{eqrbar} \end{equation}
We now have $\epsilon>0$, $\bar{r}>0$, $Y=S(\bar{x},\bar{v})$ a compact set, and $y_i\in\partial g(F(x_i))\cap Y$ with $\nabla F(x_i)^*y_i=v_i$ such that
$$\langle y_1-y_0,F(x_1)-F(x_0)\rangle\geq-\bar{r}|F(x_1)-F(x_0)|^2$$
when $|x_i-\bar{x}|<\epsilon$, $|f(x_i)-f(\bar{x})|<\epsilon$, $v_i\in\partial f(x_i)$, and $|v_i-\bar{v}|<\epsilon$.
Note that when $v_i=\nabla F(x_i)^*y_i$ then
$$\langle v_1-v_0,x_1-x_0\rangle=\langle\nabla F(x_1)^*y_i-\nabla F(x_0)^*y_1+\nabla F(x_0)^*y_1-\nabla F(x_0)^*y_0,x_1-x_0\rangle,$$ which gives
\begin{equation}
\label{eq:stat4}
\langle v_1-v_0,x_1-x_0\rangle=\langle[\nabla F(x_1)^*-\nabla F(x_0)^*]y_1,x_1-x_0\rangle+\langle\nabla F(x_0)^*(y_1-y_0),x_1-x_0\rangle.
\end{equation}
Because $Y$ is compact and $F\in\mathcal{C}^2$, we can find
	\begin{equation} \label{eqr1} r_1>0 ~\mbox{such that for all}~ y\in Y, \quad
	|[\nabla F(x_1)^*-\nabla F(x_0)^*]y|\leq r_1|x_1-x_0|.\end{equation}
Thus for all $y\in Y$,
\begin{equation}
\label{eq:stat5}
\langle[\nabla F(x_1)^*-\nabla F(x_0)^*]y,x_1-x_0\rangle\geq-r_1|x_1-x_0|^2.
\end{equation}
The final term in equation (\ref{eq:stat4}) can be written as $\langle\nabla\phi(x_0),x_1-x_0\rangle$ for the $\mathcal{C}^2$ function $\phi(x)=\langle y_1-y_0,F(x)\rangle$. Let
	\begin{equation}\label{eqr2}r_2 ~\mbox{be an upper bound for  the eigenvalues of the Hessian of the mapping}~ x\rightarrow\langle\eta,F(x)\rangle, \end{equation}
where $x$ ranges over $B_\epsilon(\bar{x})$ and $\eta$ ranges over the compact set $Y-Y$. Then in particular $\phi(x)\leq\phi(x_0)+\langle\nabla\phi(x_0),x-x_0\rangle+r_2|x-x_0|^2$ when $|x-\bar{x}|\leq\epsilon$. It follows that
\begin{equation}
\label{eq:stat6}
\begin{array}{ll}
 \langle\nabla\phi(x_0),x_1-x_0\rangle & \geq-r_2|x_1-x_0|^2+\phi(x_1)-\phi(x_0) \\
 & =-r_2|x_1-x_0|^2+\langle y_1-y_0,F(x_1)-F(x_0)\rangle\\
 & \geq-r_2|x_1-x_0|^2-\bar{r}|F(x_1)-F(x_0)|^2.\end{array}
\end{equation}
But there is also a constant $k$ such that $|F(x_1)-F(x_0)|\leq k|x_1-x_0|$ when $|x_i-\bar{x}|<\epsilon$. Using this fact with the inequalities (\ref{eq:stat5}) and (\ref{eq:stat6}) for the two terms at the end of equation (\ref{eq:stat4}) we obtain
$$\langle v_1-v_0,x_1-x_0\rangle\geq-r_1|x_1-x_0|^2-r_2|x_1-x_0|^2-\bar{r}k^2|x_1-x_0|^2.$$
Thus for
	\begin{equation}\label{eqr}
	r:=r_1+r_2+\bar{r}k^2
	\end{equation}
we have that
$$\langle v_1-v_0,x_1-x_0\rangle\geq-r|x_1-x_0|^2$$
whenever $|x_i-\bar{x}|<\epsilon$, $|f(x_i)-f(\bar{x})|<\epsilon$, $v_i\in\partial f(x_i)$, and $|v_i-\bar{v}|<\epsilon$.\qed\medskip\\
\par Now we are ready to begin with parametrized functions. Lemma \ref{lem:lemma1} demonstrates that multiplying a prox-regular function by a positive parameter does not affect its prox-regularity, a property that we shall use in the theorems that follow.
\begin{lem}
\label{lem:lemma1}
Let $f:\mathbb{R}^n\rightarrow\mathbb{R}\cup\{+\infty\}$ be prox-regular at $\bar{x}$ for $\bar{v}$ with parameters $\tilde{\epsilon}>0$ and $\tilde{r}>0$, and let $\lambda\in\mathbb{R}, \lambda>0$. Then $\lambda f(x)$ is prox-regular at $\bar{x}$ for $\lambda\bar{v}$ with parameters $\epsilon=\min\{\tilde{\epsilon},\lambda\tilde{\epsilon}\}$ and $r=\lambda\tilde{r}$.
\end{lem}
\textbf{Proof:} By prox-regularity of $f$, $\tilde{\epsilon}$ and $\tilde{r}$ are such that
$$f(x')>f(x)+\langle v,x'-x\rangle-\frac{\tilde{r}}{2}|x'-x|^2$$
whenever $x'\neq x$, $|x'-\bar{x}|<\tilde{\epsilon}$, $|x-\bar{x}|<\tilde{\epsilon}$, $|f(x)-f(\bar{x})|<\tilde{\epsilon}$, $v\in\partial f(x)$, and $|v-\bar{v}|<\tilde{\epsilon}$. Multiplying by $\lambda>0$ yields
$$\lambda f(x')>\lambda f(x)+\langle\lambda v,x'-x\rangle-\frac{\lambda\tilde{r}}{2}|x'-x|^2$$
whenever $x'\neq x$, $|x'-\bar{x}|<\tilde{\epsilon}$, $|x-\bar{x}|<\tilde{\epsilon}$, $|f(x)-f(\bar{x})|<\tilde{\epsilon}$, $v\in\partial f(x)$, and $|v-\bar{v}|<\tilde{\epsilon}$. Since $v\in\partial f(x)$, we know that $\lambda v\in\partial\lambda f(x).$ When $|v-\bar{v}|<\tilde{\epsilon}$ and $|f(x)-f(\bar{x})|<\tilde{\epsilon}$ we have $|\lambda v-\lambda\bar{v}|<\lambda\tilde{\epsilon}$ and $|\lambda f(x)-\lambda f(\bar{x})|<\lambda\tilde{\epsilon}$. Let $\epsilon=\min\{\tilde{\epsilon},\lambda\tilde{\epsilon}\}$, and let $r=\lambda\tilde{r}$. We conclude that
$$\lambda f(x')>\lambda f(x)+\langle\lambda v,x'-x\rangle-\frac{r}{2}|x'-x|^2$$
whenever $x'\neq x$, $|x'-\bar{x}|<\epsilon$, $|x-\bar{x}|<\epsilon$, $|\lambda f(x)-\lambda f(\bar{x})|<\epsilon$, $\lambda v\in\partial f(x)$, and $|\lambda v-\lambda\bar{v}|<\epsilon$.\qed\\
A direct result of Lemma \ref{lem:lemma1} is that the scalar multiplication of a prox-regular function with a positive parameter results in a para-prox-regular function.
\begin{cor}
\label{cor:cor2}
Let $f:\mathbb{R}^n\rightarrow\mathbb{R}\cup\{+\infty\}$ be prox-regular at $\bar{x}$ for $\bar{v}$ with parameters $\tilde{\epsilon}>0$ and $\tilde{r}>0$, and let $\bar{\lambda}\in\mathbb{R}$, $\bar{\lambda}>0$. Then $g(x,\lambda)=\lambda f(x)$ is para-prox-regular at $\bar{x}$ for $\bar{\lambda}\bar{v}$ with parameters $\epsilon=\min\{\tilde{\epsilon},\frac{\lambda\tilde{\epsilon}}{2}\}$ and $r=\frac{3\lambda\tilde{r}}{2}$ with compatible parametrization by $\lambda$ at $\bar{\lambda}$.
\end{cor}
\hbox{}
\textbf{Proof:} Let $\delta=\frac{\bar{\lambda}}{2}$ and consider $\lambda\in(\bar{\lambda}-\delta, \bar{\lambda}+\delta)$. Then for each such $\lambda$, using $\epsilon_\lambda=\min\{\tilde{\epsilon},\delta\tilde{\epsilon}\}$ we have that
$$g(x',\lambda)>g(x,\lambda)+\langle\lambda v,x'-x\rangle-\frac{\lambda\tilde{r}}{2}|x'-x|^2$$
whenever $x'\neq x$, $|x'-\bar{x}|<\epsilon_\lambda$, $|x-\bar{x}|<\epsilon_\lambda$, $|g(x,\bar{\lambda})-g(\bar{x},\bar{\lambda})|<\epsilon_\lambda$, $\lambda v\in\partial_xg(x,\lambda)$, and $|\lambda v-\bar{\lambda}\bar{v}|<\epsilon_\lambda$. Let $\epsilon=\min\{\delta,\tilde{\epsilon},\delta\tilde{\epsilon}\}$ and $r=\max\{\lambda\tilde{r}\}$. Then we have
$$g(x',\lambda)>g(x,\lambda)+\langle\lambda v,x'-x\rangle-\frac{r}{2}|x'-x|^2$$
when $x'\neq x$, $|x'-\bar{x}|<\epsilon$, $|x-\bar{x}|<\epsilon$, $|g(x,\lambda)-g(\bar{x},\bar{\lambda})|<\epsilon$, $\lambda v\in\partial_xg(x,\lambda)$, $|\lambda v-\bar{\lambda} v|>\epsilon$, and $|\lambda-\bar{\lambda}|<\epsilon$.\qed\medskip\\
Next we have a minor extension of \cite[Theorem 3.2]{calcprox}, where Poliquin and Rockafellar showed that, under certain conditions, the sum of two prox-regular functions is again prox-regular. We extend this to the sum of any finite number of prox-regular functions. The result also extracts the prox-regularity parameters from the proof, thus allowing their use in later results.
\begin{thm}
\label{thm:main1}
For $i\in\{1,2,\ldots,m\}$, let $f_i:\mathbb{R}^n\rightarrow\mathbb{R}\cup\{+\infty\}$. Let $\bar{x}\in\bigcap\limits_{i=1}^m\dom f_i$ and assume that the only choice of $v_i\in\partial f^\infty_i(\bar{x})$ with $\sum\limits_{i=1}^mv_i=0$ is $v_i=0$ for all $i=1,2,\ldots,m$. Define $f(x):=\sum\limits_{i=1}^mf_i(x)$ and let $\bar{v}\in\partial f(\bar{x})$. If $f_i$ is prox-regular at $\bar{x}$ for $v_i$ with parameters $\epsilon_i>0$ and $r_i>0$, for each $v_i\in\partial f_i(\bar{x})$ such that $\sum\limits_{i=1}^mv_i=\bar{v}$, then $f$ is prox-regular at $\bar{x}$ for $\bar{v}$ with parameters $\epsilon=\min_i\{\epsilon_i\}$ and $r=m\max_i\{r_i\}$.
\end{thm}
\textbf{Proof:} Define $g(u_1,u_2,\ldots,u_m):=\sum\limits_{i=1}^mf_i(u_i)$ and $F(x):=(x,x,\ldots,x)$. Then $g$ is lsc and proper, $F\in\mathcal{C}^2$, and $f(x)=g(F(x))$. Notice that $\nabla F(x)=(I,I,\ldots,I)^T$, so the only possible $y\in\partial^\infty g(F(\bar{x}))$ with $\nabla F(\bar{x})^*y=0$ is $y=0$. Therefore the constraint qualification on $y$ in Theorem \ref{thm:prelim} holds. We have that $g$ is prox-regular at $F(\bar{x})$ for every $y\in\partial g(F(\bar{x}))$ with $\nabla F(\bar{x})^*y=\bar{v}$, so all the conditions of Theorem \ref{thm:prelim} are satisfied. Let $\epsilon$ and $r$ be as described in the proof of that theorem. In particular, $r_1=0$ in this case since $\nabla F(x_1)-\nabla F(x_0)=0$. Also $r_2=0$, since $\langle\eta,F(x)\rangle=\langle\eta,(x,x,\ldots,x)\rangle=x\eta$, which implies $\nabla\langle\eta,F(x)\rangle=\eta$ and $\nabla^2\langle\eta,F(x)\rangle=0$.  Finally,
$$|F(x_1)-F(x_0)|=|(x_1-x_0,x_1-x_0,\ldots,x_1-x_0)|=\sqrt{m(x_1-x_0)^2}=\sqrt{m}|x_1-x_0|,$$ which gives $\lambda=\sqrt{n}$. Hence $r=m\bar{r}=m\max\limits_i\{r_i\}$. We conclude that $f$ is prox-regular at $\bar{x}$ for $\bar{v}$ with these parameters.\qed
\begin{cor}
\label{cor:cor1}
For $i\in\{1,2,\ldots,m\}$, let $f_i:\mathbb{R}^n\rightarrow\mathbb{R}\cup\{+\infty\}$ be prox-regular at $\bar{x}$ for $v_i\in\partial f_i(\bar{x})$ with parameters $\epsilon_i>0$ and $r_i>0$. Let $\bar{\lambda}=(\lambda_1,\lambda_2,\ldots,\lambda_m)>0$ (i.e. $\lambda_i>0$ for all $i$), and define $f(x):=\sum\limits_{i=1}^m\lambda_if_i(x).$ Then $f$ is prox-regular at $\bar{x}$ for $\bar{\lambda}\bar{v}$ with parameters $\epsilon=\min\limits_i\{\epsilon_i,\lambda_i\epsilon_i\}$ and $r=m\max\limits_i\{\lambda_ir_i\}$, where $\bar{v}=\sum\limits_{i=1}^m\lambda_iv_i\in\partial_xf(\bar{x},\bar{\lambda})$.
\end{cor}

\textbf{Proof:} Apply Lemma \ref{lem:lemma1} to each $f_i$ in Theorem \ref{thm:main1}.\qed\medskip\\
We now examine parametrized functions by multiplying each of the prox-regular functions $f_i$ in Theorem \ref{thm:main1} by a parameter $\lambda_i$, and we show that the weighted sum of prox-regular functions is a para-prox-regular function.

\begin{thm}
\label{thm:par}
For $i\in\{1,2,\ldots,m\}$, let $f_i:\mathbb{R}^n\rightarrow\mathbb{R}\cup\{+\infty\}$. Let $\bar{\lambda}=(\bar{\lambda}_1,\bar{\lambda}_2,\ldots,\bar{\lambda}_m)>0$. Let $\bar{x}\in\bigcap\limits_{i=1}^m\dom f_i$ and assume that
\begin{equation}
\label{eq:cq}
\mbox{the only choice of }v_i\in\partial f^\infty_i(\bar{x})\mbox{ with }\sum\limits_{i=1}^mv_i=0\mbox{ is }v_i=0\mbox{ for all }i.
\end{equation}
Define
$$f(x,\lambda):=\sum\limits_{i=1}^m\lambda_if_i(x)$$
where $\lambda=(\lambda_1,\lambda_2,\ldots,\lambda_m)$, and let $\bar{v}\in\partial_xf(\bar{x},\bar{\lambda})$. If $f_i$ is prox-regular at $\bar{x}$ for all $\bar{v}_i\in\partial f_i(\bar{x})$ such that $\sum\limits_{i=1}^m\bar{\lambda}_i\bar{v}_i=\bar{v}$, with parameters $\epsilon_i$ and $r_i$, then $f$ is para-prox-regular at $\bar{x}$ for $\bar{v}$ with compatible parametrization by $\lambda$ at $\bar{\lambda},$ with parameters $\epsilon=\min\limits_i\{\frac{\bar{\lambda}_i}{2},\epsilon_i,\frac{\bar{\lambda}_i}{2}\epsilon_i\}$ and $r=m\max\limits_i\{\frac{3\bar{\lambda}_i}{2}r_i\}.$
\end{thm}
\textbf{Proof:} Since $\bar{\lambda}_i>0$ and $f_i$ is prox-regular at $\bar{x}$ for $\bar{v}_i$ for all $i$, Lemma \ref{lem:lemma1} applies and we have that each $\lambda_i f_i$ is prox-regular at $\bar{x}$ for $\bar{\lambda}_i\bar{v}_i$ with parameters $\tilde{\epsilon}_i=\min\{\epsilon_i,\bar{\lambda}_i\epsilon_i\}$ and $\tilde{r}_i=\bar{\lambda}_i r_i$. Since $\bar{\lambda}_i>0$, we have $\partial^\infty\bar{\lambda}_if_i=\partial^\infty f_i$ \cite[equation 10(6)]{rockwets}, so condition (\ref{eq:cq}) holds for $v_i\in\partial^\infty\bar{\lambda}_if_i$. Thus by Corollary \ref{cor:cor1}, we have that $f(x,\bar{\lambda})$ is prox-regular in $x$ at $\bar{x}$ for $\bar{v}$ with parameters $\bar{\epsilon}=\min\limits_i\{\tilde{\epsilon}_i\}$ and $\bar{r}=m\max\limits_i\{\tilde{r}_i\}.$ (Note that here we are viewing $f$ as a function of $x$, with $\bar{\lambda}$ fixed.) Define $\delta:=\frac{1}{2}\min\limits_i\{\bar{\lambda}_i\}$, and $\mathcal{B}:=\{\lambda:\lambda_i\in[\bar{\lambda}_i-\delta,\bar{\lambda}_i+\delta]\mbox{ for all }i\}.$ Since $\lambda>0$ for all $\lambda\in\mathcal{B}$, by similar arguments as above we have that $f(x,\lambda)$ is prox-regular in $x$ at $\bar{x}$ for $\bar{v}$ for every fixed $\lambda\in\mathcal{B}$, with parameters $\epsilon_\lambda=\min\limits_i\{\epsilon_i, \lambda_i\epsilon_i\}$ and $r_\lambda=m\max\limits_i\{\lambda_i r_i\}.$ But $\min\limits_{\lambda\in\mathcal{B}}\{\lambda_i\}=\delta$ and $\max\limits_{\lambda\in\mathcal{B}}\{\lambda_i\}=\max\limits_i\{\bar{\lambda}_i\}+\delta,$ so for all $\lambda\in\mathcal{B}$ we have that $f(x,\lambda)$ is prox-regular in $x$ at $\bar{x}$ for $\bar{v}$ with parameters $\hat{\epsilon}=\min\{\epsilon_i,\delta\epsilon_i\}$ and $\hat{r}=m\max\limits_i\{(\bar{\lambda}_i+\delta)r_i\}.$ Since $\bar{\lambda}_i+\delta\leq\frac{3\bar{\lambda}_i}{2},$ we may overestimate and set $r=m\max\limits_i\{\frac{3\bar{\lambda}_i}{2}r_i\}.$ Therefore we may conclude that $f$ is para-prox-regular at $(\bar{x},\bar{\lambda})$ for $\bar{v}$ with parameters $\epsilon=\min\{\delta,\hat{\epsilon}\}=\min\limits_i\{\frac{\bar{\lambda}_i}{2},\epsilon_i,\frac{\bar{\lambda}_i}{2}\epsilon_i\}$ and $r=m\max\limits_i\{\frac{3\bar{\lambda}_i}{2}r_i\}.$\qed
\begin{ex}
Let $f_1:\mathbb{R}^n\rightarrow\mathbb{R}\cup\{\infty\}$ and $f_2:\mathbb{R}^n\rightarrow\mathbb{R}\cup\{\infty\}$ be prox-regular at $\bar{x}\in\bigcap\limits_{i=1}^2\dom f_i.$ Assume the only choice of $\hat{v}_1\in\partial f_1^\infty(\bar{x})$ and $\hat{v}_1\in\partial f_2^\infty(\bar{x})$ with $\hat{v}_1+\hat{v}_2=0$ is $\hat{v}_1=\hat{v}_2=0.$ Define
$$f(x,\lambda):=\lambda f_1(x)+(1-\lambda)f_2(x).$$
Let $\bar{v}\in\partial_xf(\bar{x},\bar{\lambda})$ such that $\bar{v}=\lambda v_1+(1-\lambda)v_2$ with $v_1\in\partial f_1(\bar{x})$ and $v_2\in\partial f_2(\bar{x})$. Then $f$ is para-prox-regular at $\bar{x}$ for $\bar{v}=\lambda v_1+(1-\lambda)v_2$ with compatible parametrization by $\lambda$ at any $\bar{\lambda}\in(0,1).$
\end{ex}
\textbf{Details:} For $f_1$ and $f_2$ lower-$\mathcal{C}^2$ functions, this result has already been proved \cite{proxave}. The lower-$\mathcal{C}^2$ requirement can be generalized to prox-regular by using Theorem \ref{thm:par}.\qed\medskip\\
The following proposition and lemma are a restatement of a result from \cite{stability}. We shall use them in the proof of Theorem \ref{thm:max}, para-prox-regularity of a weighted finite max of $\mathcal{C}^0$ functions.
\begin{prop} \cite[Prop. 2.2]{stability}
\label{prop:sa}
Suppose that $f:\mathbb{R}^n\times\mathbb{R}^m$ is strongly amenable in $x$ at $\bar{x}$ with compatible parametrization by $\lambda$ at $\bar{\lambda}$, in the sense that on some neighborhood of $(\bar{x},\bar{\lambda})$ there is a composite representation $f(x,\lambda)=g(F(x,\lambda))$ in which $F:\mathbb{R}^n\times\mathbb{R}^d\rightarrow\mathbb{R}^m$ is a $\mathcal{C}^2$ mapping and $g:\mathbb{R}^m\rightarrow\mathbb{R}$ is a convex, proper, lsc function for which $F(\bar{x},\bar{\lambda})\in D:=\dom g$ and
$$z\in N_D(F(\bar{x},\bar{\lambda})), \nabla_xF(\bar{x},\bar{\lambda})^*z=0\Rightarrow z=0.$$
Then as long as $\bar{v}\in\partial_xf(\bar{x},\bar{\lambda})$, one has $f$ continuously para-prox-regular in $x$ at $\bar{x}$ for $\bar{v}$ with compatible parametrization by $\lambda$ at $\bar{\lambda}$.
\end{prop}
\begin{cor}
\label{cor:cor100}
For $i\in\{1,2,\ldots,m\}$, let $f_i:\mathbb{R}^n\rightarrow\mathbb{R}$ be $\mathcal{C}^2.$ Let $\bar{\lambda}=(\bar{\lambda}_1,\bar{\lambda}_2,\ldots,\bar{\lambda}_m)>0$. Let $\bar{x}\in\bigcap\limits_{i=1}^m\dom f_i$ and $\bar{v}\in\partial_xf(\bar{x},\bar{\lambda})$. Define
$$f(x,\lambda):=\max\limits_i\{\lambda_if_i(x)\}$$
where $\lambda=(\lambda_1,\lambda_2,\ldots,\lambda_m)$. Then $f$ is continuously para-prox-regular in $x$ at $\bar{x}$ for $\bar{v}$ with compatible parametrization by $\lambda$ at $\bar{\lambda}$.
\end{cor}
\textbf{Proof:} Let
$$F(x,\lambda)=(\lambda_1f_1(x),\lambda_2f_2(x),\ldots,\lambda_mf_m(x))$$
and
$$g(u_1,u_2,\ldots,u_m)=\max\{u_1,u_2,\ldots,u_m\}.$$
Then $F$ is $\mathcal{C}^2$, since each component is a product of two $\mathcal{C}^2$ functions, $g$ is proper lsc convex, and $f(x,\lambda)=g(F(x))$.\\
Therefore $f$ is strongly amenable. Since $D=\dom g$ is $\mathbb{R}^n$, $N_D(F(\bar{x},\bar{\lambda}))=\{0\}$, so the constraint qualification of Proposition \ref{prop:sa} holds. We have all the requirements of Proposition \ref{prop:sa}, and its conclusion is our desired result.\qed\medskip\\
Our next result is a relaxation of the $\mathcal{C}^2$ condition of Corollary \ref{cor:cor100}. It is sufficient that the $f_i$ be $\mathcal{C}^1$, prox-regular functions. In relaxing the result, we provide a direct proof and also extract the prox-regularity parameters from the proof.

\begin{thm}
\label{thm:max}
For $i\in\{1,2,\ldots,m\}$, let $f_i:\mathbb{R}^n\rightarrow\mathbb{R}$ be $\mathcal{C}^1$. Let $\bar{\lambda}=(\bar{\lambda}_1,\bar{\lambda}_2,\ldots,\bar{\lambda}_m)>0$. Let $\bar{x}\in\bigcap\limits_{i=1}^m\dom f_i$. Define
$$f(x,\lambda):=\max\limits_i\{\lambda_if_i(x)\}$$
where $\lambda=(\lambda_1,\lambda_2\ldots,\lambda_m)$, and let $\bar{v}\in\partial_xf(\bar{x},\bar{\lambda})$. If the $f_i$ are prox-regular at $\bar{x}$ for each $v_i\in\partial f_i(\bar{x})$ with parameters $\epsilon_i>0$ and $r_i>0$, then $f$ is para-prox-regular at $\bar{x}$ for $\bar{v}$ with parameters $\epsilon=\min\limits_i\{\epsilon_i,\frac{\bar{\lambda}_i\epsilon_i}{2}\}$ and $r=\max\limits_i\{\frac{3\bar{\lambda}_ir_i}{2}\}$, with compatible parametrization by $\lambda$ at $\bar{\lambda}$.
\end{thm}
\textbf{Proof:} Let $A(\bar{x},\bar{\lambda})=\{i:f(\bar{x},\bar{\lambda})=\bar{\lambda}_if_i(\bar{x})\}$ be the active set of indices. Then $\partial_xf(\bar{x},\bar{\lambda})=\conv\limits_{i\in A(\bar{x},\bar{\lambda})}\{\bar{\lambda}_i\nabla f_i(\bar{x})\}$ \cite[Theorem 2.8.2]{optanal}, so every $\bar{v}\in\partial_xf(\bar{x},\bar{\lambda})$ has the form $\sum\limits_{i\in A(\bar{x},\bar{\lambda})}\theta_i\bar{\lambda}_iv_i$ for some $v_i\in\partial f_i(\bar{x})$ and some set of non-negative $\{\theta_i\}_{i\in A(\bar{x},\bar{\lambda})}$ such that $\sum\limits_{i\in A(\bar{x},\bar{\lambda})}\theta_i=1$. Corollary \ref{cor:cor2} tells us that the $\lambda_if_i$ are para-prox-regular at $(\bar{x},\bar{\lambda}_i)$ for all $\bar{\lambda}_iv_i\in\partial\bar{\lambda}_if_i(\bar{x})$. In particular, for each $i$ we can take $\hat{\epsilon}_i=\min\{\epsilon_i, \frac{\bar{\lambda}_i\epsilon_i}{2}\}$ and $\hat{r}_i=\frac{3\bar{\lambda}_ir_i}{2}$ to get
\begin{equation}
\label{equation:max9}
\lambda_if_i(x')\geq\lambda_if_i(x)+\langle\lambda_i\hat{v}_i,x'-x\rangle-\frac{\hat{r}_i}{2}|x'-x|^2
\end{equation}
for $x'\neq x$, $|x'-\bar{x}|<\hat{\epsilon}_i$, $|x-\bar{x}|<\hat{\epsilon}_i$, $|\lambda_if_i(x)-\bar{\lambda}_if(\bar{x})|<\hat{\epsilon}_i$, $\lambda_i\hat{v}_i\in\partial\lambda_if_i(x)$, $|\lambda_i\hat{v}_i-\bar{\lambda}_iv_i|<\hat{\epsilon}_i$, and $|\lambda_i-\bar{\lambda}_i|<\hat{\epsilon}_i.$ This is a system of $m$ inequalities, so using $\epsilon=\min\limits_i\{\hat{\epsilon}_i\}$ and $r=\max\limits_i\{\hat{r}_i\}$ we can say that for all $i\in\{1,2,\ldots,m\}$
\begin{equation}
\label{equation:max10}
\lambda_if_i(x')\geq \lambda_if_i(x)+\langle\lambda_i\hat{v}_i,x'-x\rangle-\frac{r}{2}|x'-x|^2
\end{equation}
for all $x'\neq x, |x'-\bar{x}|<\epsilon, |x-\bar{x}|<\epsilon, |\lambda_i-\bar{\lambda}_i|<\epsilon, \lambda_i\hat{v}_i\in\partial\lambda_if_i(x), |\lambda_i\hat{v}_i-\bar{\lambda}_iv_i|<\epsilon$, $|\lambda_if_i(x)-\bar{\lambda}_if_i(\bar{x})|<\epsilon$ (with $r$ and $\epsilon$ no longer dependent on $i$). Since $f(x',\lambda)\geq\lambda_if_i(x')$ for all $i$, we can replace the left-hand side of inequality (\ref{equation:max10}) with $f(x',\lambda)$ to get
\begin{equation}
\label{equation:max11}
f(x',\lambda)\geq \lambda_if_i(x)+\langle\lambda_i\hat{v}_i,x'-x\rangle-\frac{r}{2}|x'-x|^2
\end{equation}
for all $x'\neq x, |x'-\bar{x}|<\epsilon, |x-\bar{x}|<\epsilon, |\lambda_i-\bar{\lambda}_i|<\epsilon, \lambda_i\hat{v}_i\in\partial\lambda_if_i(x), |\lambda_i\hat{v}_i-\bar{\lambda}_iv_i|<\epsilon$, $|\lambda_if_i(x)-\bar{\lambda}_if_i(\bar{x})|<\epsilon.$
Let $x',x,\lambda,v$ be such that $x'\neq x$, $|x'-\bar{x}|<\epsilon$, $|x-\bar{x}|<\epsilon$, $v\in\partial_xf(x,\lambda)$, $|v-\bar{v}|<\epsilon$, $|\lambda-\bar{\lambda}|<\epsilon$, $|f(x,\lambda)-f(\bar{x},\bar{\lambda})|<\epsilon$. Let $A(x,\lambda)=\{i:f(x,\lambda)=\lambda_if_i(x)\}$. Then $\partial_xf(x,\lambda)=\conv\limits_{i\in A(x,\lambda)}\{\lambda_i\partial f_i(x)\}$. So $v\in\partial_xf(x,\lambda)$ has the form $\sum\limits_{i\in A(x,\lambda)}\theta_i\lambda_i\hat{v}_i$ for some $\hat{v}_i\in\partial f_i(x)$ and some non-negative $\{\theta_i\}_{i\in A(x,\lambda)}$ such that $\sum\limits_{i\in A(x,\lambda)}\theta_i=1$. We can also substitute $f(x,\lambda)$ on the right-hand side of inequality (\ref{equation:max11}) for those $i\in A(x,\lambda)$ to get
\begin{equation}
\label{equation:max12}
f(x',\lambda)\geq f(x,\lambda)+\langle\lambda_i\hat{v}_i,x'-x\rangle-\frac{r}{2}|x'-x|^2.
\end{equation}
\smallskip
Multiplying (\ref{equation:max12}) by $\theta_i$ and adding over $i\in A(x,\lambda)$ yields
$$\sum\limits_{i\in A(x,\lambda)}\theta_if(x',\lambda)\geq\sum\limits_{i\in A(x,\lambda)}\theta_if(x,\lambda)+\langle\sum\limits_{i\in A(x,\lambda)}\theta_i\lambda_i\hat{v}_i,x'-x\rangle-\sum\limits_{i\in A(x,\lambda)}\theta_i\frac{r}{2}|x'-x|^2$$
which simplifies to
\begin{equation}
\label{equation:max13}
f(x',\lambda)\geq f(x,\lambda)+\langle v,x'-x\rangle-\frac{r}{2}|x'-x|^2.
\end{equation} Therefore $f$ is para-prox-regular at $\bar{x}$ for $\bar{v}$ with compatible parametrization by $\lambda$ at $\bar{\lambda}$, with parameters $\epsilon$ and $r$.\qed

\section{Open Questions}
\label{sec:open}
There are two areas of further research that are made apparent by this work. First, our results concerning weighted parametrized sums and weighted parametrized max functions rely on starting with a finite number of prox-regular functions. A natural direction to continue with the development of these theorems would be to explore the case where the family of functions is infinite. Second, all the results found in this paper are for functions in finite-dimensional space. A natural extension of this research would be to reexamine these results in a Hilbert or Banach space.

\subsection*{Acknowledgement}

The authors thank the anonymous referee and our colleagues, H.\ Bauschke, P.\ Loewen, and X.\ Wang, for constructive feedback during the final stages of revision of this paper.

\end{document}